\documentclass[11pt]{article}
\usepackage{graphicx, colordvi}
\usepackage{dsfont}
\usepackage{epsfig}
\usepackage{mathrsfs}
\usepackage{amssymb,amsfonts,amsmath,amsthm,cite,color}
\usepackage{dsfont}
\usepackage{epsfig}
\usepackage{mathrsfs}
\usepackage{longtable}
\usepackage{enumerate}
\usepackage{geometry}
\usepackage{hyperref}
\hypersetup{
	colorlinks=true,
	linkcolor=black,
	filecolor=black,
	urlcolor=black,
	citecolor=black,
}

\newtheorem{theo}{Theorem}

\newtheorem{coro}[theo]{Corollary}

\newtheorem{lem}[theo]{Lemma}

\makeatletter \@addtoreset{equation}{section}
\@addtoreset{theo}{section} \makeatother

\newcommand{\bN} { {\mathbb{N}}}

\newcommand{\bQ} { {\mathbb{Q}}}
\newcommand{\bZ} { {\mathbb{Z}}}

\newcommand{\bK} { {\mathbb{K}}}

\DeclareMathOperator{\ann}{ann}
\DeclareMathOperator{\ord}{ord}

\def\qed{\hfill \rule{4pt}{7pt}}
\def\pf{\noindent {\it Proof.} }

\author{}
\title{}
\begin{document}
\begin{center}

{\large \bf Power-partible Reduction and Congruences for Schr\"oder Polynomials}
\end{center}
\begin{center}
{Chen-Bo Jia}$^{1}$, {Rong-Hua Wang}$^{2}$, {Michael X.X. Zhong}$^{3}$

  $^{1,2}$School of Mathematical Sciences\\
   Tiangong University \\
   Tianjin 300387, P.R. China \\
   $^1$2230111379@tiangong.edu.cn\\
   $^2$wangronghua@tiangong.edu.cn
   \\[10pt]

  $^3$School of Science\\
   Tianjin University of Technology \\
   Tianjin 300384, P.R. China\\
   zhong.m@tjut.edu.cn
\end{center}

\vskip 6mm \noindent {\bf Abstract.}
In this note, we apply the power-partible reduction to show the following arithmetic properties of large Schr\"oder polynomials $S_n(z)$ and little Schr\"oder polynomials $s_n(z)$:
for any odd prime $p$, nonnegative integer $r\in\bN$, $\varepsilon\in\{-1,1\}$ and $z\in\bZ$ with $\gcd(p,z(z+1))=1$, we have
\[
\sum_{k=0}^{p-1}(2k+1)^{2r+1}\varepsilon^k S_k(z)\equiv 1\pmod {p}\quad \text{and}
\quad
\sum_{k=0}^{p-1}(2k+1)^{2r+1}\varepsilon^k s_k(z)\equiv 0\pmod {p}.
\]	

\section{Introduction}

For $n\in\bN$, the $n$th large Schr\"oder number $S_n$ and little Schr\"oder number $s_n$ are defined respectively by
\[
S_n=\sum_{k=0}^{n}\binom{n}{k}\binom{n+k}{k}\frac{1}{k+1}
\quad
\text{and}
\quad
s_n=\sum_{k=1}^{n}\frac{1}{n}\binom{n}{k}\binom{n}{k-1}2^{n-k}.
\]
It is well-known that $S_n=2s_n$ for $n\in\bN^{*}$.
Note also that $S_0=1$ and $s_0=0$.

There are many combinatorial objects that are counted by $S_n$ or $s_n$.
One can consult items \href{https://oeis.org/A006318}{A006318} and \href{https://oeis.org/A001003}{A001003} in the \href{https://oeis.org/}{OEIS} \cite{OEIS} for some of their classical combinatorial interpretations.

Arithmetic properties of Schr\"oder numbers are also extensively investigated.
For example,
Sun \cite{Sun2011} proved \begin{equation*}
\sum_{k=1}^{p-1}\frac{S_k}{m^k}\equiv
\frac{m^2-6m+1}{2m}\left(1-
\left(\frac{m^2-6m+1}{p}\right)\right) \pmod p,
\end{equation*}
where $p$ is an odd prime, $(-)$ is the Legendre symbol and $m$ is any integer not divisible by $p$.
Cao and Pan \cite{CaoPan2017} showed that, for $n\geq 1$ and $\alpha \geq 1$,
\[
S_{n+2^\alpha}\equiv S_n+2^{\alpha+1} \pmod {2^{\alpha+2}}.
\]
In 2016, Liu \cite{Liu2016} confirmed the supercongruence
\[
\sum_{k=1}^{p-1}D_kS_k\equiv 2p^3B_{p-3}-2pH^{*}_{p-1} \pmod {p^4},
\]
originally conjectured by Sun \cite{Sun2011}.
Here $D_n$ denotes the Delannoy number, $B_n$ denotes the Bernoulli number and $H^{*}_n=
\sum_{k=1}^{n}\frac{(-1)^k}{k}$
 is the alternating harmonic number.

Motivated by the definition of Schr\"oder numbers, Sun \cite{Sun2012, Sun2018} defined the large Schr\"oder polynomials $S_n(z)$ and the little Schr\"oder polynomials $s_n(z)$ as
\begin{equation*}
S_n(z)=\sum_{k=0}^{n}\binom{n}{k}\binom{n+k}{k}\frac{1}{k+1}z^k
\quad \text{and} \quad s_n(z)=\sum_{k=1}^{n}\frac{1}{n}\binom{n}{k}\binom{n}{k-1}z^{k-1}(z+1)^{n-k}.
\end{equation*}
Clearly, $S_n(1)$ and $s_n(1)$ reduce to the $n$th large and little Schr\"oder numbers respectively.
Sun \cite[Lemma 2.1]{Sun2018} proved that
\begin{equation}\label{eq:s_n(z)}
(z+1)s_n(z)=S_n(z),\quad\forall n\in\bN^*
\end{equation}
and
\begin{equation}\label{eq:D_n}
2z(2n+1)S_n(z)=D_{n+1}(z)-D_{n-1}(z),\quad\forall n\in\bN^*,
\end{equation}
where $D_n(z)$ are the $n$th central Delannoy polynomials given by
\[
D_n(z)=\sum_{k=0}^{n}\binom{n}{k}\binom{n+k}{k}z^k.
\]

In recent years, Sun \cite{Sun2014,Sun2016,Sun2018} presented many amazing congruences involving $D_n(z)$, $s_n(z)$ or Ap\'ery numbers $A_n$.
For example, Sun \cite{Sun2018} derived that
\[
\sum_{k=0}^{p-1}kD_k(z)s_{k+1}(z)\equiv
2(z(z+1))^{(p-1)/2} \pmod p.
\]
In 2014,
Sun \cite{Sun2014} proved that
\[
\sum_{k=0}^{p-1}(2k+1)D_k(z)\equiv \frac{p}{z}((z+1)^p-1) \pmod {p^3},
\]
if $z\in\bZ$ and $p$ is an odd prime with $\gcd(p,z)=1$.
Recently, Xia and Sun \cite{XiaSun2023} confirmed the conjecture posed by Sun \cite{Sun2016} that
for each $r\in\bN$ and prime $p>3$, there is a $p$-adic integer $c_r$ only depending on $r$ such that
\[
\sum_{k=0}^{p-1}(2k+1)^{2r+1}(-1)^kA_k\equiv c_r p \left(\frac{p}{3}\right) \pmod {p^3}.
\]
For positive integers $m$ and $\alpha$, the arithmetic property of
\[
\sum_{k=0}^{n-1}\varepsilon^k(2k+1)A_k^{(\alpha)}(z)^m
\]
has also been discussed by Guo and Zeng \cite{GuoZeng2012} and Pan \cite{Pan2014}.
Here $\varepsilon\in\{-1,1\}$ and $A_k^{(\alpha)}(x)$ is the generalized Ap\'ery polynomial, which is also called Schmidt polynomials in \cite{GuoZeng2012}.

In 2021, Hou, Mu and Zeilberger \cite{HouMuZeil2021} introduced an algorithmic process, which is called the power-partible reduction, to derive congruences for hypergeometric terms mechanically.
Recently, Wang and Zhong \cite{WZ2023+} generalized the reduction to holonomic sequences and utilized it to study the arithmetic properties of
\[\sum_{k=0}^{p-1}(2k+1)^{2r}\varepsilon^k F_k
\quad\text{and}\quad
\sum_{k=0}^{p-1}(2k+1)^{2r+1}\varepsilon^kF_k,\]
where $\varepsilon\in\{-1,1\}$ and $F_k$ are Ap\'ery numbers $A_k$ or central Delannoy polynomials $D_k(z)$.

By Zeilberger's algorithm \cite{Zeil1990,Zeilberger1990b, Zeilberger1991} and the definition of power-partibility given in \cite{WZ2023+}, it is straightforward to check that both $S_n(z)$ and $s_n(z)$ are power-partible.
However, the original technique used in \cite{WZ2023+} loses its power when we consider the arithmetic properties of $\sum\limits_{k=0}^{p-1}(2k+1)^{2r+1}\varepsilon^k S_k(z)$ or $\sum\limits_{k=0}^{p-1}(2k+1)^{2r+1}\varepsilon^k s_k(z)$.

In this note, we extend the technique of the power-partible reduction and then use it to deduce the following congruences for $S_n(z)$ and $s_n(z)$.
\begin{theo}\label{th:(2k+1)^{2r+1}Sk(z)}
Let $p$ be an odd prime, $\varepsilon\in\{-1,1\}$ and $z\in\bZ$ with $\gcd(p,z(z+1))=1$.
Then for all $r\in\bN$, we have
\begin{equation}\label{1q:(2k+1)^{2r+1}}
\sum_{k=0}^{p-1}(2k+1)^{2r+1}\varepsilon^k S_k(z)\equiv 1\pmod {p}
\ \text{ and } \
\sum_{k=0}^{p-1}(2k+1)^{2r+1}\varepsilon^k s_k(z)\equiv 0\pmod {p}.
\end{equation} 	
\end{theo}

Taking $z=1$ in Theorem \ref{th:(2k+1)^{2r+1}Sk(z)} yields the following congruences for Schr\"oder numbers.

\begin{coro}\label{coro:Sk(1) 3.5}
Let $p$ be an odd prime and $\varepsilon\in\{-1,1\}$. Then for all $r\in\bN$, we have
\begin{equation*}
\sum_{k=0}^{p-1}(2k+1)^{2r+1}\varepsilon^kS_k\equiv 1\pmod {p}
\quad \text{and}\quad
\sum_{k=0}^{p-1}(2k+1)^{2r+1}\varepsilon^ks_k\equiv 0\pmod {p}.
\end{equation*} 		
\end{coro}

The rest of the note is organized as follows.
Section 2 sums up basic techniques of power-partible reduction developed in \cite{WZ2023+} and provides an extended version (Theorem 2.1).
A proof of Theorem \ref{th:(2k+1)^{2r+1}Sk(z)} based on the extended reduction technique is presented in Section 3.

\section{Extended power-partible reduction}
Let $\bK$ be a field of characteristic $0$ and $\bK[k]$ the polynomial ring over $\bK$.
The set of \emph{annihilators} of $F(k)$ is defined by
\begin{equation*}\label{eq:Annihilator}
\ann F(k):=\left\{L=\sum_{i=0}^{J}a_i(k)\sigma^i\in \bK[k][\sigma]\mid L(F(k))=0\right\},
\end{equation*}
where $J\in\bN=\{0,1,2,\ldots\}$ and $\sigma$ is the shift operator (that is, $\sigma F(k)=F(k+1)$).

A sequence $F(k)$ is said to be \emph{holonomic} (or, \emph{P-recursive}) if and only if $\ann F(k)\neq \{0\}$.
The \emph{order} of $L=\sum_{i=0}^{J}a_i(k)\sigma^i\in \bK[k][\sigma]$ is defined to be $J$ if $a_J(k)\neq 0$, denoted by $\ord(L)=J$.
The minimum order of all $L\in\ann F(k)\setminus\{0\}$ is called the \emph{order} of $F(k)$.
A holonomic sequence $F(k)$ of order $J$ is called \emph{summable} if
\begin{equation*}
F(k)=\Delta\left(\sum_{i=0}^{J-1}u_i(k)F(k+i)\right)
\end{equation*}
holds for some rational functions $u_i(k)\in\bK(k)$, where $\Delta=\sigma-1$ is the difference operator.

Given an operator $L=\sum\limits_{i=0}^{J}a_i(k)\sigma^i\in \bK[k][\sigma]$ of order $J$,
the \emph{adjoint} of $L$ is defined by
$
L^{\ast}=\sum_{i=0}^{J}\sigma^{-i}a_i(k),
$
that is,
\[
L^{\ast}(x(k))=\sum_{i=0}^{J}a_i(k-i)x(k-i)
\]
for any $x(k)\in\bK[k]$.
In 2018, van der Hoeven \cite{Hoeven2018} showed that,
if $L=\sum_{i=0}^{J}a_i(k)\sigma^i\in\ann F(k)$,
then for any $x(k)\in\bK[k]$ one has
\begin{equation}\label{eq: summable}	L^{\ast}(x(k))F(k)=\Delta\left(-\sum_{i=0}^{J-1}u_i(k)F(k+i)\right),
\end{equation}
where
\begin{equation*}
u_i(k)=\sum_{j=1}^{J-i}a_{i+j}(k-j)x(k-j),\quad i=0,1,2,\ldots,J-1.
\end{equation*}
Summing over $k$ from $0$ to $n-1$ on both sides of identity \eqref{eq: summable} leads to
\begin{equation}\label{eq:rec ann}
\sum_{k=0}^{n-1}L^{\ast}(x(k)) F(k)
=\left(\sum_{i=0}^{J-1}u_i(0)F(i)\right)
-\left(\sum_{i=0}^{J-1}u_i(n)F(n+i)\right).
\end{equation}

The \emph{degree} of $L=\sum\limits_{i=0}^{J}a_i(k)\sigma^i\in \bK[k][\sigma]$, denoted by $\deg L$, is given by
\begin{equation*}\label{eq:d and bk}
\deg L=\max_{0\leq \ell \leq J} \{\deg b_\ell(k)-\ell\},
\end{equation*}
where $b_{\ell}(k)=\sum_{j={\ell}}^{J}\binom{j}{\ell}a_{J-j}(k+j-J)$.
Let
\begin{equation*}\label{eq:nonnegative roots}
R_{L}=\{s\in\bN \mid \sum_{\ell=0}^{J}[k^{d+\ell}](b_\ell(k))s^{\underline{\ell}}=0\}.
\end{equation*}
Here $[k^{d+\ell}](b_\ell(k))$ denotes the coefficient of $k^{d+\ell}$ in $b_\ell(k)$ and $s^{\underline{\ell}}$ denotes the falling factorial $s(s-1)\cdots(s-\ell+1)$.
Then $L$ is called \emph{nondegenerated} if $R_{L}=\emptyset$, and \emph{degenerated} otherwise.

When $L=\sum\limits_{i=0}^{J}a_i(k)\sigma^i\in \bK[k][\sigma]$ is nondegenerated and there exists a $\gamma\in\bK$ such that
\begin{equation*}\label{eq:power-reduce-condition}
a_i(\gamma+k)=(-1)^{\deg L} a_{J-i}(\gamma-k-J),\quad i=0,1,\ldots,\lfloor\frac{J}{2}\rfloor,
\end{equation*}
then $L$ is called \emph{power-partible} with respect to $\gamma$ \cite{WZ2023+}.
If $L\in\ann F(k)$, we also say $F(k)$ is power-partible with respect to $\gamma$.

Wang and Zhong \cite[Lemma 2.5]{WZ2023} showed that
\begin{equation*}\label{eq:degree-L*}
\deg L^{\ast}(x(k))=\deg x(k)+\deg L,
\quad\forall x(k)\in\bK[k],
\end{equation*}
if $L$ is nondegenerated.
If further $L$ is power-partible with respect to $\gamma$, then we can show that $L^*(x(k))$ is a linear combination of purely odd powers $(k-\gamma)^{2i+1}$ or purely even powers $(k-\gamma)^{2j}$ for those $x(k)\in\bK[k]$ satisfying certain symmetry condition.

\begin{theo}\label{th:reduction}
Suppose $L=\sum_{i=0}^{J}a_i(k)\sigma^i\in\bK[k][\sigma]\setminus\{0\}$ is of order $J$ and degree $d$, and power-partible with respect to some $\gamma\in\bK$.
Let $\ell\in\bN$ and $x_j(k)\in\bK[k]$ ($j=\ell,\ell+1,\ldots$) be a series of polynomials with $\deg x_j(k)=j$ such that
\begin{equation}\label{eq:symmetry}
x_j(\gamma+k)=(-1)^{\deg x_j(k)} x_j(\gamma-k-J).
\end{equation}
Then $L^*(x_j(k))$ can be written as a linear combination of $(k-\gamma)^t$ over $\bK$ with $t\equiv d+j\pmod 2$.
Therefore, for any $m\in\bN^*$, there exist some $u_i,v_j\in\bK$ such that
\begin{equation}\label{eq:power-partible reduction}
(k-\gamma)^m=\sum_{\substack{0\le i<d+\ell\\ i\equiv m \pmod 2}}u_i(k-\gamma)^i
+\sum_{\substack{\ell\le j\le m-d\\ d+j\equiv m \pmod 2}}v_jL^\ast(x_{j}(k)).
\end{equation}
\end{theo}
\pf Theorem \ref{th:reduction}, a slight generalization of Theorem 2.4 and 2.6 in \cite{WZ2023+}, follows directly from the proof of Theorem 2.4 in \cite{WZ2023+}.
\qed

We call \eqref{eq:power-partible reduction} a \emph{power-partible reduction} of $(k-\gamma)^m$ using $L$ and $x_j(k)$, $j\geq\ell$.
Next, we will discuss the structure of those polynomials $x(k)$ satisfying the symmetry condition \eqref{eq:symmetry}.
Let
\[
x_1(k)=k-\gamma+\frac{J}{2}
\quad\text{and}\quad
x_2(k)=(k-\gamma+\frac{J}{2})^2-q
\]
with $q\in\bK$.
Then it is clear that \eqref{eq:symmetry} holds for both $x_1(k)$ and $x_2(k)$.
Notice that the symmetry condition \eqref{eq:symmetry} is multiplicative, namely, if \eqref{eq:symmetry} holds for $x(k)$ and $y(k)$, then it must also hold for the product $x(k)y(k)$.
Therefore, the $x_j(k)$ satisfying \eqref{eq:symmetry} can be chosen as the product
\begin{equation}\label{eq:form}
\alpha\prod_{i}\tilde{x}_i(k)
\end{equation}
where $\alpha\in\bK$ is a constant and $\tilde{x}_i(k)$ equals $x_1(k)$ or $x_2(k)$.

The following Lemma shows that if we allow $q$ to be an element in the algebraic closure $\overline{\bK}$ of $\bK$, then all polynomials satisfying \eqref{eq:symmetry} are of the form \eqref{eq:form}.

\begin{lem}\label{lem: symmetry}
Let $\gamma\in\bK$, $J\in\bN^*$, $x(k)\in\bK[k]$. Then we have
\begin{equation}\label{eq:symmetry2}
x(\gamma+k)=(-1)^{\deg x(k)} x(\gamma-k-J)
\end{equation}
if and only if $x(k)$ can be written as a constant multiple of a product of finitely many polynomials of
\begin{equation}\label{eq: basic form}
k-\gamma+\frac{J}{2}
\quad\text{and}\quad
(k-\gamma+\frac{J}{2})^2-q
,
\end{equation}
where $q\in\overline{\bK}$.
\end{lem}
\pf The sufficiency is clear from the discussion above.
We proceed to show the necessity.
Suppose $x(k)\in\bK[k]$ satisfies \eqref{eq:symmetry2}.
If $x(k)$ is a constant, the conclusion is trivially true.

Now suppose $x(k)$ is not a constant and $\tilde{k}\in\overline{\bK}$ is one of its roots.
Then we can write $\tilde{k}=\gamma-\frac{J}{2}-\tilde{q}$ for some $\tilde{q}\in\overline{\bK}$.
By \eqref{eq:symmetry2}, one can see $\gamma-\frac{J}{2}+\tilde{q}$ must also be a root of $x(k)$.
If $\tilde{q}=0$, then $k-\gamma+\frac{J}{2}$ is a factor of $x(k)$;
otherwise,
\begin{equation}\label{eq:q}
(k-\gamma+\frac{J}{2}+\tilde{q})(k-\gamma+\frac{J}{2}-\tilde{q})
=(k-\gamma+\frac{J}{2})^2-\tilde{q}^2
\end{equation}
is a factor of $x(k)$.
In either case, we can divide $x(k)$ by the factor to get a new polynomial still satisfying \eqref{eq:symmetry2} but with a lower degree.
Taking $q=\tilde{q}^2$ in \eqref{eq:q}, the necessity then follows by induction.
\qed

When doing the power-partible reduction in \cite{WZ2023+}, $x_j(k)$ were chosen to be constant multiples of $(k-\gamma+\frac{J}{2})^j$.
The appearance of $x_2(k)=(k-\gamma+\frac{J}{2})^2-q$ here allows more flexibility, which finally leads to a proof of Theorem \ref{th:(2k+1)^{2r+1}Sk(z)} in Section 3.

\section{New congruences for Schr\"oder polynomials}

Recall that the large Schr\"oder polynomials $S_n(z)$ and the little Schr\"oder polynomials $s_n(z)$ are defined respectively by
\begin{equation*}
S_n(z)=\sum_{k=0}^{n}\binom{n}{k}\binom{n+k}{k}\frac{1}{k+1}z^k\text{ and } s_n(z)=\sum_{k=1}^{n}\frac{1}{n}\binom{n}{k}\binom{n}{k-1}z^{k-1}(z+1)^{n-k}.
\end{equation*}
In this section, $z$ is taken to be a parameter in $\bZ$.
We will assume $z(z+1)\neq 0$ since $p$ is an odd prime and $\gcd(p,z(z+1))=1$ in Theorem.
Then by the connection \eqref{eq:s_n(z)}, we know $S_n(z)$ and $s_n(z)$ have the same recurrence relations and thus $\ann S_n(z)=\ann s_n(z)$ if we take $\bK=\bQ$ to be the field of rational numbers.

To simplify the discussion, we will take $\varepsilon\in\{-1,1\}$ and
\[
F_n(z)=\varepsilon^n S_n(z)\text{ or }\varepsilon^n s_n(z).
\]
By Zeilberger's algorithm, we find that
\begin{equation}\label{eq: S L}
L=a_2(k)\sigma^2+a_1(k)\sigma+a_0(k)\in \ann F_k(z),
\end{equation}
where
\begin{equation*}
a_2(k)=k+3, a_1(k)=-\varepsilon(2k+3)(1+2z), a_0(k)=k.
\end{equation*}
It is easy to check that $J=\ord L=2$, $d=\deg L=1$, $L$ is nondegenerated and
\begin{equation*}
a_0(\gamma+k)=(-1)^d a_2(\gamma-k-J)
 \quad\text{and}\quad
a_1(\gamma+k)=(-1)^d a_1(\gamma-k-J),
\end{equation*}
for $\gamma=-\frac{1}{2}$.
That is, $F_k(z)$ is power-partible with respect to $\gamma=-\frac{1}{2}$.
By the identity \eqref{eq:rec ann}, for any polynomial $x_s(k)\in\bK[k]$ of degree $s$, we have
\begin{equation}\label{eq:x(k)Sk(z)}
\sum_{k=0}^{n-1}L^{\ast}(x_s(k))F_k(z)
=(u_0(0)F_0(z)+u_1(0)F_1(z))
-\left(u_0(n)F_n(z)+u_1(n)F_{n+1}(z)\right),
\end{equation}
where $u_0(n)=(n+1)x_s(n-2)-\varepsilon(2n+1)(1+2z)x_s(n-1)$ and $u_1(n)=(n+2)x_s(n-1).$

To obtain congruences after reduction later, we need to eliminate the initial term $u_0(0)F_0(z)+u_1(0)F_1(z)=\left(x(-2)
-\varepsilon(1+2z)x_s(-1)\right)F_0(z)+2x_s(-1)F_1(z)$.
Since $\gamma=-\frac{1}{2}$ and $J=2$.
Taking $q=1/4$ in \eqref{eq: basic form}, we have $(k+\frac{3}{2})^2-\frac{1}{4}=(k+1)(k+2)$ is one of the basic forms in \eqref{eq: basic form}.
Then by Lemma \ref{lem: symmetry} we can choose
\begin{equation}\label{eq:y_s+2}
x_{s+2}(k)=2(2k+3)^s(k+1)(k+2),\quad s\in\bN
\end{equation}
to proceed the power-partible reduction in \eqref{eq:power-partible reduction} and eliminate the initial term in \eqref{eq:x(k)Sk(z)} at the same time.
Note that $d=1$ and $\ell=2$.
By Theorem \ref{th:reduction}, for each $r\in\bN$, there are constants $c_s, c\in\bK$ such that
\begin{equation*}\label{eq:reT(2k+1)}
(2k+1)^{2r+1}=\sum_{s=0}^{r-1}c_sL^\ast(x_{2s+2}(k))+c(2k+1).
\end{equation*}
The following lemma gives a more concrete description of the constants $c_s$ and $c$.
\begin{lem}\label{lem:(2k+1)^{2r+1}}
Let $L$ be given by \eqref{eq: S L} and $x_{s+2}(k)$ given by \eqref{eq:y_s+2}. Then for any $r\in\bN$ we have
\begin{equation}\label{eq:(2k+1)^{2r+1}}
(2k+1)^{2r+1}=\sum_{s=0}^{r-1}
\frac{v_s}{\eta^{u_s}}L^\ast(x_{2s+2}(k))+(2k+1),
\end{equation}
for some $u_s\in\bN$ and $v_s\in\bZ$, where $\eta=\frac{1-\varepsilon(1+2z)}{2}\in\bZ$. 	
\end{lem}
\pf
For simplicity, let $\ell=2k+1$. By the definition of $L^\ast$, we have
\begin{align}\label{eq:ylt}
&L^\ast(x_{s+2}(k))=\sum_{i=0}^{2}a_{i}(k-i)x_{s+2}(k-i)\notag\\
=&2k(k+1)\left((k+2)(2k+3)^s-\varepsilon(2k+1)(1+2z)(2k+1)^s+(k-1)(2k-1)^s\right)\\
=&\frac{1}{4}(\ell-1)(\ell^2+4\ell+3)(\ell+2)^s+\frac{1}{4}(\ell+1)(\ell^2-4\ell+3)(\ell-2)^s-\frac{\varepsilon(1+2z)(\ell^2-1)\ell^{s+1}}{2}\notag\\
=&(\ell^3-\ell)\sum_{\substack{j=0\\ j\text{ even}}}^{s}\binom{s}{j}2^{j-1}\ell^{s-j}+(3\ell^2-3)\sum_{\substack{j=1\\j\text{ odd}}}^{s}\binom{s}{j}2^{j-1}\ell^{s-j}-\frac{\varepsilon(1+2z)(\ell^2-1)\ell^{s+1}}{2}\notag\\
=&\eta\ell^{s+3}-\eta\ell^{s+1}+\sum_{\substack{j=2\\j\text{ even}}}^{s}\binom{s}{j}2^{j-1}\ell^{s-j}(\ell^3-\ell)+(3\ell^2-3)\sum_{\substack{j=1\\j\text{ odd}}}^{s}\binom{s}{j}2^{j-1}\ell^{s-j}\notag\\
=&\eta\ell^{s+3}+\sum_{j=1}^{\lfloor \frac{s+3}{2}\rfloor}e_j\ell^{s+3-2j}, \nonumber
\end{align}
where $\eta=\frac{1-\varepsilon(1+2z)}{2}\in\bZ$ and $e_j\in\bZ$ for all $j=1,2,\ldots,\lfloor \frac{s+3}{2}\rfloor$.
By the expression for $L^\ast(x_{s+2})$, the power-partible reduction on $(2k+1)^{2r+1}$ reveals that
\begin{equation}\label{eq:(2k+1)^{2r+1}decomp}
(2k+1)^{2r+1}=\sum_{s=0}^{r-1}
\frac{v_s}{\eta^{u_s}}L^\ast(x_{2s+2}(k))+\frac{v}{\eta^{u}}(2k+1),
\end{equation}
for some $u_s,u\in\bN$, $v_s,v\in\bZ$.

By the expression \eqref{eq:ylt} for $L^\ast(x_{s+2}(k))$, we know $k$ divides $L^\ast(x_{s+2}(k))$, namely, $k=0$ is a root of $L^\ast(x_{s+2})$.
Letting $k=0$ in \eqref{eq:(2k+1)^{2r+1}decomp}, one can see $\frac{v}{\eta^{u}}=1$.
This completes the proof.
\qed

Now we consider the arithmetic properties of $\sum_{k=0}^{p-1} L^{\ast}(x_{2s+2}(k))F_k(z)$ and $\sum_{k=0}^{p-1}(2k+1)F_k(z)$, respectively.

\begin{lem}
Let $L$ be given by \eqref{eq: S L}.
Then for any $n,z\in\bZ$ with $n>1$ and polynomial $y(k)\in\bZ[k]$, we have
\begin{equation*}
\sum_{k=0}^{n-1} L^{\ast}(y(k)(k+1)(k+2))F_k(z) \equiv 0 \pmod {n(n^2-1)}.
\end{equation*}
\end{lem}
\pf Let $x(k)=y(k)(k+1)(k+2)$, equality \eqref{eq:x(k)Sk(z)} reduces to
\begin{equation}\label{eq:L*y(k)}
\sum_{k=0}^{n-1} L^{\ast}(y(k)(k+1)(k+2)) F_k(z)
=-\left(u_0(n)F_n(z)+u_1(n)F_{n+1}(z)\right),
\end{equation}
here $u_0(n)=n(n+1)((n-1)y(n-2)-\varepsilon(2n+1)(1+2z)y(n-1))$ and $u_1(n)=n(n+1)(n+2)y(n-1).$
Since $L\in\ann F_k(z)=0$, apparently
\begin{equation}\label{eq:L(S_n(z))}
(n+2)F_{n+1}(z)=\varepsilon(2n+1)(1+2z)F_n(z)-(n-1)F_{n-1}(z),\quad\forall n\ge 1.
\end{equation}
Substituting identity \eqref{eq:L(S_n(z))} into identity \eqref{eq:L*y(k)}, we derive
\begin{equation}\label{eq:L*(x(k))Sk(z) v2}
	\sum_{k=0}^{n-1} L^{\ast}(y(k)(k+1)(k+2)) F_k(z)
	=n(n^2-1)\left(y(n-1)F_{n-1}(z)-y(n-2)F_{n}(z)\right).
\end{equation}
This completes the proof.
\qed

Letting $y(k)=2(2k+3)^s$ in \eqref{eq:L*(x(k))Sk(z) v2}, we get the following corollary.
\begin{coro}\label{lem:L*(x(k))Sk(z) tong}
Let $L$ be given by \eqref{eq: S L} and $x_s(k)$ given by \eqref{eq:y_s+2}.
Then for any $z\in\bZ$ and $n,s\in\bN$ with $n>1$, we have
\begin{equation*}\label{eq:L*(x(k))Sk(z) tong}
\sum_{k=0}^{n-1} L^{\ast}(x_{s+2}(k))F_k(z) \equiv 0 \pmod {2n(n^2-1)}.
\end{equation*}
\end{coro}
The next lemma is the last piece of the puzzle.
\begin{lem}\label{lem:3.2}
Let $p$ be an odd prime, $z\in\bZ$ with $\gcd(p,z(z+1))=1$.
Then we have
\begin{equation*}
\sum_{k=0}^{p-1}(2k+1)\varepsilon^k S_k(z)\equiv 1 \pmod {p} \quad\text{and}\quad
\sum_{k=0}^{p-1}(2k+1)\varepsilon^k s_k(z)\equiv 0 \pmod {p},
\end{equation*}	
where $\varepsilon\in\{-1,1\}$.
\end{lem}
\pf
By equality \eqref{eq:D_n}, we have
\begin{align*}
2z\sum_{k=0}^{p-1}(2k+1)\varepsilon^k S_k(z)&=2z+\sum_{k=1}^{p-1}\varepsilon^k2z(2k+1)S_k(z)\\
&=2z+\sum_{k=1}^{p-1}\varepsilon^k(D_{k+1}(z)-D_{k-1}(z))\\
&=2z+D_p(z)+\varepsilon D_{p-1}(z)-D_1(z)-\varepsilon D_0(z).\\
\end{align*}
Since  $D_0(z)=1,D_1(z)=1+2z$, we deduce that
\begin{equation}\label{eq: 2(2k+1)Sk}
2z\sum_{k=0}^{p-1}(2k+1)\varepsilon^kS_k(z)=D_p(z)+\varepsilon D_{p-1}(z)-1-\varepsilon.	
\end{equation}
It is clear to see that
\begin{align*}
D_p(z)&=\sum_{k=0}^{p}\binom{p}{k}\binom{p+k}{k}z^k\\
&=1+\sum_{k=1}^{p-1}\binom{p}{k}\binom{p+k}{k}z^k+\binom{2p}{p}z^p\\
&\equiv 1+\binom{2p}{p}z^p \pmod {p}\\
&\equiv 1+2z \pmod {p},
\end{align*}
and
\begin{align*}
D_{p-1}(z)=\sum_{k=0}^{p-1}\binom{p-1}{k}\binom{p-1+k}{k}z^k
=1+\sum_{k=1}^{p-1}\binom{p-1}{k}\binom{p-1+k}{k}z^k
\equiv 1 \pmod{p}.
\end{align*}
From the above two congruences and identity \eqref{eq: 2(2k+1)Sk}, we obtain
\begin{equation}\label{eq: initial conruence}
2z\sum_{k=0}^{p-1}(2k+1)\varepsilon^kS_k(z)\equiv 2z\pmod{p}.
\end{equation}
Since $p$ is an odd prime and $\gcd(p,z)=1$, the first congruence in Lemma \ref{lem:3.2}  follows from \eqref{eq: initial conruence}.
By the fact that $S_k(z)=(1+z)s_k(z)$ for all $k\in\bN^*$ and $S_0(z)=1,s_0(z)=0$, it follows from \eqref{eq: initial conruence} that
\[2z(1+z)\sum_{k=0}^{p-1}(2k+1)\varepsilon^k s_k(z)\equiv 0\pmod{p}.\]
The second congruence then follows from the condition that $p$ is an odd prime and $\gcd(p,z(z+1))=1$.
\qed

\noindent\emph{Proof of Theorem \ref{th:(2k+1)^{2r+1}Sk(z)}}:
Multiplying both sides of identity \eqref{eq:(2k+1)^{2r+1}} with $F_k(z)$ and then
summing over $k$ from $0$ to $p-1$, we get
\begin{equation*}\sum_{k=0}^{p-1}(2k+1)^{2r+1}F_k(z)\equiv \sum_{k=0}^{p-1}(2k+1)F_k(z) \pmod {p}
\end{equation*}
with the help of Corollary \ref{lem:L*(x(k))Sk(z) tong} and the fact that $\gcd(p,\eta)=1$.
Congruence \eqref{1q:(2k+1)^{2r+1}} then follows immediately from Lemma \ref{lem:3.2}.
\qed

\noindent \textbf{Acknowledgments.}
This work was supported by the National Natural Science Foundation of China (No. 12101449, 12271511, 12271403) and the Natural Science Foundation of Tianjin, China (No. 22JCQNJC00440).

\end{document}